\documentclass[openright,a4paper,american,amsart,11pt]{article}

\usepackage[latin1]{inputenc}

\usepackage{amsmath}
\usepackage[english]{babel}

\usepackage{amsfonts}
\input xy
\xyoption{all}
\usepackage[dvips]{graphicx}
\usepackage{color}

\newtheorem{thm}{Theorem}[section]

\newtheorem{prop}[thm]{Proposition}
\newtheorem{lemma}[thm]{Lemma}

\newcommand{\findemo}{\hfill
                     $\Box$ \vspace{1.5 ex}}
\newcommand{\RR}{\mathbb R}
\newcommand{\CC}{\mathbb C}

\newcommand{\cL}{\mathcal L}

\begin{document}

\selectlanguage{english}

\title{A NON-ALGEBRAIC PATCHWORK}
\author{Benoît Bertrand \and  Erwan Brugallé}
\date{}
\maketitle
\begin{abstract}
Itenberg and Shustin's pseudoholomorphic curve patchworking is in
principle more flexible than Viro's original algebraic one. It was
natural to wonder if the former method allows one to construct
non-algebraic objects.  In this paper we construct the first examples
of patchworked real pseudoholomorphic curves in $\Sigma_n$ whose
position with respect to the pencil of lines cannot be realised by
any homologous real algebraic curve.
\end{abstract}

\setlength{\parindent}{0mm}

\section{Introduction}

Viro's patchworking has been since the seventies one of the most important and
fruitful methods in topology of real algebraic varieties. It was
applied in the proof of a lot of meaningful results in this field
(e.g.  \cite{V1}, \cite{V4}, \cite{Sh2}, \cite{Lop},\cite{I3},
\cite{I1}, \cite{Haa},\cite{Br2}, \cite{B1}, \cite{IV1}, \cite{Ber1},
\cite{Mik1}, and  \cite{Sh3}). 
Here we will only consider the case of curves in 
$\CC P^2$ and in rational geometrically ruled surfaces $\Sigma_n$.

Viro Method allows to construct an algebraic curve $A$ out of simpler
 curves $A_i$ so that the topology of $A$ can be deduced from the
 topology of the initial curves $A_i$. Namely one gets a curve with
 Newton polytope $\Delta$ out of curves whose Newton polytopes are the
 $2$-simplices of a subdivision $\sigma$ of $\Delta$, and one can see
 the curve $A$ as a gluing of the curves $A_i$. Moreover, if all the
 curves $A_i$ are real, so is the curve $A$.  For a detailed account
 on the Viro Method we refer to \cite{V1} and \cite{V4} for example.

One of the hypotheses of the Viro Method is that $\sigma$ should be
convex (i.e.  the $2$-simplices of $\sigma$ are the domains of
linearity of a piecewise linear convex function).  In the original
Patchworking Theorem, one also requires the curves $A_i$ to be totally
nondegenerate. In particular, one can only glue nonsingular
curves. Later, E. Shustin proved in \cite{Sh1} that, under some
numerical conditions depending on the types and number of the
singularities, it is possible to patchwork singular curves keeping the
singular points.

On the other hand, I. Itenberg and E. Shustin proved in \cite{IS} a
pseudoholomorphic patchworking theorem:
 they showed that applying the Viro Patchworking with any subdivision
(non necessarily convex) and with reduced curves $A_i$ with arbitrary
singularities, one can glue the $A_i$'s, keeping singular points, to
obtain a \textit{real pseudoholomorphic curve}. More precisely, given
some (maybe singular) curves $A_i$ whose Newton polygons are the
$2$-simplices of a subdivision of the quadrangle with vertices
$(0,0)$, $(k,0)$, $(k,l)$ and $(l+nk,0)$, Itenberg and Shustin gave a
way to construct a pseudoholomorphic curve $C$ of bidegree $(k,l)$ in
the rational geometrically ruled surface $\Sigma_n$, whose position
with respect to the pencil of lines can be deduced from the initial
curves $A_i$.
Isotopy types  realizable by (algebraic or
pseudoholomorphic) curves obtained via a patchworking procedure are
called \textit{patchworked curves}.

\vspace{2ex} Pseudoholomorphic curves were introduced by M. Gromov in
\cite{Gro} to study symplectic $4$-mani\-folds.  A real
pseudoholomorphic curve $C$ on $\mathbb CP^2$ or $\Sigma_n$ is an
immersed Riemann surface which is a $J$-holomorphic curve in some tame
almost complex structure $J$ such that the exceptional section (in
$\Sigma_n$ with $n\ge 1$) is $J$-holomorphic, $conj(C)=C$, and
$conj_*\circ J_p=-J_p\circ conj_*$ (where $conj$ is the standard
complex conjugation and $p$ is any point of $C$).  It as been realized
since then that real pseudoholomorphic curves share a lot of
properties with real algebraic ones (see for example \cite{OS1} and
\cite{OS2}). It is still unknown if there exist nonsingular real
pseudoholomorphic curves in $\CC P^2$ or in $\Sigma_n$ which are
isotopic to no homologous real algebraic curves (this is the so called
real symplectic isotopy problem).  Note that, not requiring the
exceptional section in $\Sigma_n$ to be $J$-holomorphic,
J-Y. Welschinger constructed in \cite{W} examples of real
pseudoholomorphic curves on $\Sigma_n$ for $n\ge 2$ which are not
isotopic to any real algebraic curve realizing the same homology
class.

\vspace{2ex} In the surfaces $\mathbb R \Sigma_n$, there is a natural
pencil of lines $\mathcal L$, and one can study curves there up to
fiberwise isotopy. Two curves $C_1$ and $C_2$ in $\mathbb R \Sigma_n$
are said to be \textit{$\mathcal L$-isotopic} if there exists an
isotopy $\phi(t,x)$ of $\mathbb R \Sigma_n$ mapping $C_1$ to $C_2$
such that for any $t\in[0;1]$, for any $p\in C_1$ and for any fiber
$F$ of $\mathbb R \Sigma_n$, $\phi(t,F)$ is a fiber of $\mathbb R
\Sigma_n$, and the  intersection multiplicity of $C_1$ and $F$ at $p$ is the
 intersection multiplicity of $\phi(t,C_1)$ and $\phi(t,F)$ at $\phi(t,p)$.

A lot of examples are known of nonsingular real pseudoholomorphic
curves in $\mathbb R \Sigma_n$ which are $\mathcal L$-isotopic to no
homologous real algebraic curves (see for example \cite{OS1},
\cite{OS2}, \cite{Br1}). However, as far as we know none of those
examples are constructed with the pseudoholomorphic patchworking of
Itenberg and Shustin, and the question of the existence of a
patchworked pseudoholomorphic curves with any kind of non-algebraic
behaviour was open.

\vspace{2ex} In \cite{Br3} we proved that in the case of curves of
bidegree $(3,0)$ in $\Sigma_n$, the patchworked pseudoholomorphic
curve is always isotopic to a real algebraic one in the same homology
class.

In this paper we construct the first examples of patchworked real
pseudoholomorphic curves in $\Sigma_n$ whose position with respect to
the pencil of lines cannot be realised by any homologous real
algebraic curve.

\vspace{2ex}
A smooth curve $C$ in $\mathbb  R
\Sigma_n$ is said to be \textit{$\mathcal L$-nonsingular} if $C$
intersects any fiber transversally, except for a 
finite number of fibers which have an ordinary  tangency point with
one of the branches of $C$, and intersect transversally the other
branches of $C$. 

A smooth curve $C$ in $\mathbb  R
\Sigma_n$ which is $\mathcal L$-singular is called \textit{smooth
  $\mathcal L$-singular}.

\begin{thm}
For any $d \ge 3$ there exists a smooth
  $\mathcal L$-singular real  pseudoholomorphic
patchworked curve of bidegree $(d,0)$ in $\Sigma_2$ which is not
  $\mathcal L$-isotopic to any  real algebraic curve in $\Sigma_2$ of the
same bidegree.
\end{thm}

The pseudoholomorphic construction is done in
Proposition \ref{pseudo}, and the algebraic obstruction is proved in
Proposition \ref{alg ob}.

\section{Rational geometrically ruled surfaces}\label{defi}

In this section we fix our notations for the surfaces $\Sigma_n$.
\textit{The $n^{th}$ rational geometrically ruled surface}, denoted by
$\Sigma_n$, is the surface obtained by taking four copies of $\mathbb
C^2$ with coordinates $(x,y)$, $(x_2,y_2)$, $(x_3,y_3)$ and
$(x_4,y_4)$, and by gluing them along $(\mathbb C^*)^2$ with the
identifications $(x_2,y_2)=(1/x,y/x^n)$, $(x_3,y_3)=(x,1/y)$ and
$(x_4,y_4)=(1/x,x^n/y)$. Let us denote by $E$ (resp. $B$ and $F$) the
algebraic curve in $\Sigma_n$ defined by the equation $\{y_3=0\}$
(resp. $\{y=0\}$ and $\{x=0\}$). The coordinate system $(x,y)$ is
called \textit{standard}.  The projection $\pi$~: $(x,y)\mapsto x$ on
$\Sigma_n$ defines a $\mathbb CP^1$-bundle over $\mathbb CP^1$.  The
intersection numbers of $B$ and $F$ are respectively $B\circ B=n$,
$F\circ F=0$ and $B\circ F=1$. The surface $\Sigma_n$ has a natural
real structure induced by the complex conjugation in $\mathbb C^2$,
and the real part $\RR\Sigma_n$ of $\Sigma_n$ is a torus if $n$ is
even and a Klein bottle if $n$ is odd. The restriction of $\pi$ on
$\mathbb R\Sigma_n$ defines a pencil of lines denoted by $\mathcal L$.

 The group $H_2(\Sigma_n,\mathbb Z)$ is isomorphic to $\mathbb
 Z\times\mathbb Z$ and is generated by the classes of $B$ and
 $F$. Moreover, one has $E=B-nF$. An algebraic   or pseudoholomorphic
 curve on $\Sigma_n$ is  
said to be of \textit{bidegree} $(k,l)$ if it realizes the homology
 class $kB+lF$ in $H_2(\Sigma_n,\mathbb Z)$. Its equation in
 $\Sigma_n\setminus E$ is
\begin{displaymath}
\sum_{i=0}^k a_{k-i}(X,Z)Y^i \eqno{(*)}
\end{displaymath}
where $a_j(X,Z)$ is a homogeneous polynomial of degree $nj+l$.

\section{Proof of the main Theorem}\label{patchwork}

Our constructions use the Pseudoholomorphic Patchworking Theorem (see
\cite{IS}). The fact that the curves constructed have a position with
respect to $\cL$ which is not realizable by an algebraic curve comes
from a condition on the degree of a certain univariate polynomial.

\begin{prop}\label{perturb tgce}
Let $d\ge 3$ be a natural number and note $k=\left[ \frac{d}{2}\right]$. Choose
$k$ real numbers $0<\alpha_1<\alpha_2<\ldots < \alpha_k$ and define the curve 
 $C_{sing}$ by the equation
$$\begin{array}{ll}
\prod_{i=1}^k \left(Y^2 - \alpha_i X  \right) & \textrm{if d is
  even}
\\
\\ Y\prod_{i=1}^k \left(Y^2 - \alpha_i X  \right) & \textrm{if d is
  odd}.
\end{array}$$
Then, for $\epsilon>0$ small enough, the  real algebraic curve
$C=C_{sing}+\epsilon 
(Y^d + X   - X^{d-1})$ in $\CC P^2$ satisfies
\begin{itemize}
\item $C$ is nonsingular,
\item at the point  $[0:0:1]$ (resp. $[1:0:0]$), the curve $C$ is
  locally given by the equation  $Y^d+X=0$ (resp. $Y^d-X^{d-1}=0$),
\item there exists a real number $a>0$ such that the line of equation
  $X-aZ=0$ intersects the curve $C$ in $d$ distinct real points.
\end{itemize}  
\end{prop}
\textit{Proof. } Straightforward.\findemo

Let us denote by $\widetilde C$ the curve defined by the equation $X^d
C(X,\frac{Y}{X})$. Then, the Newton polygon of $\widetilde C$ is the
triangle with vertices $(0,d)$, $(d+1,0)$ and $(2d-1,0)$ and the curve
$\widetilde C$ is locally given by the equation $Y^d+ X^{d+1}=0$
(resp. $Y^d-X^{2d-1}=0$) at the point $[0:0:1]$
(resp. $[1:0:0]$). Moreover, there exists a real number $b>0$ such
that the line of equation $X-bZ=0$ intersects the curve $\widetilde C$
in $d$ distinct real points.

\vspace{2ex} Define the polynomials $P_1=Y^d+X +1$,
$P_2=Y^d-X^{d-1}+X^d$, $P_3= Y^d + X^d + X^{d+1}$ and $P_4= Y^d-
X^{2d-1}+X^{2d} $. The curves $P_1$ and $P_3$ (resp. $P_2$ and $P_4$)
have a maximal tangency point with the line $X+Z=0$ (resp. $X-Z=0$).

\begin{figure}[h]
      \centering
 \begin{tabular}{c}
 \includegraphics[width=6cm, angle=0]{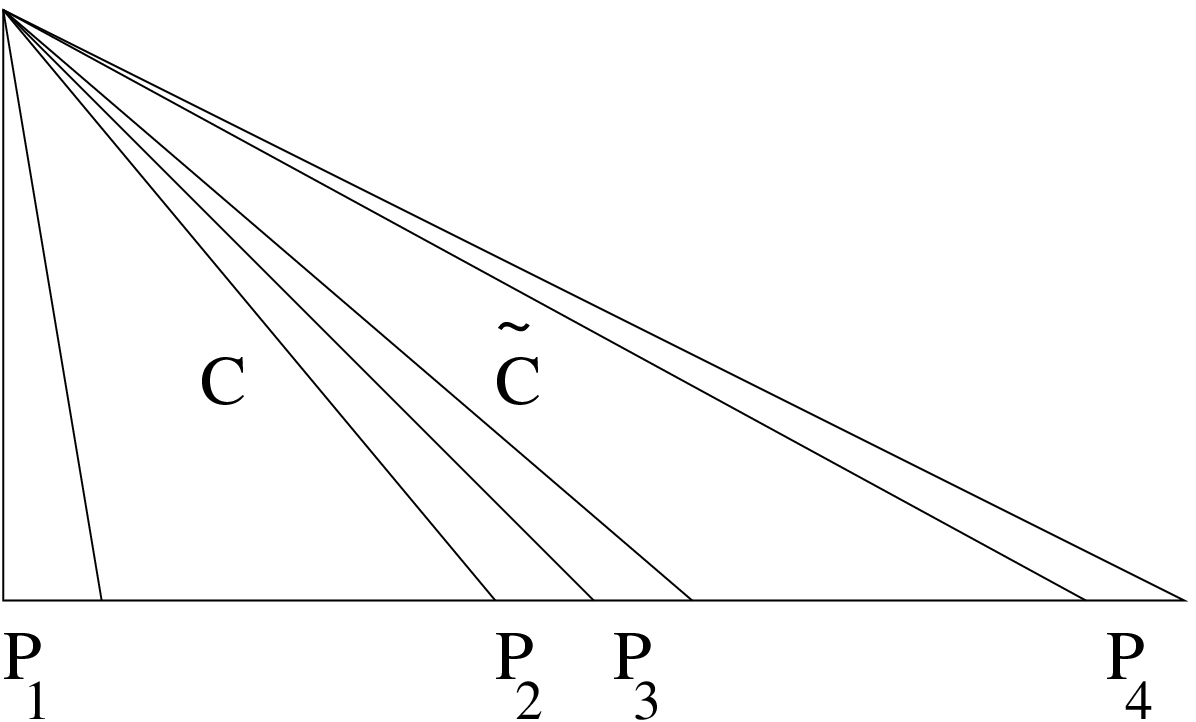}
 \end{tabular}
\caption{}
 \label{patch}
\end{figure}

\vspace{2ex}
One can now patchwork the polynomials $C$,  $\widetilde C$, $P_1$,
$P_2$, $P_3$, and  $P_4$ (see Figure \ref{patch}).
According to Itenberg and Shustin's theorem (see \cite{IS}),
 one can glue pseudoholomorphically all the pieces of the patchwork
keeping the tangency conditions with respect to $\cL$. This proves the
following proposition.
\begin{prop}\label{pseudo}
There exists a patchworked nonsingular real pseudoholomorphic curve
$\Gamma$ of bidegree $(d,0)$ in $\Sigma_2$ such that there exist real
numbers $x_1<x_2<0<z_1<x_4<z_2<x_6$ satisfying
\begin{itemize}
\item $\Gamma$ has a maximal tangency point with the fibers $X-x_iZ=0$
with $i=1,\ldots 4$,
\item $\Gamma$ has $d$ distinct real intersection points with the
  fibers $X-z_iZ=0$ with $i=1,2$.\findemo
\end{itemize}  
\end{prop}

\begin{lemma}\label{a2}
Let $P(X,Y)$ be real a polynomial of the form $$Y^d + a_{2}(X)Y^{d-2}
+ a_{3}(X)Y^{d-3} + \ldots + a_{d}(X).$$ and $x$ be a real number.  If
the univariate polynomial $Q_x(Y)=P(x,Y)$ has only real roots then
$a_2(x) \le 0$ and $a_2(x)= 0$ if and only if $0$ is a root of order
$d$.
\end{lemma}

\textit{proof.}  Let $y_i$ be the $d$ roots of $Q_x$.  The sum
$\sum_{i=1}^d y_i$ is $0$ and the second coefficient satisfies
$a_{2}=\sum_{i<j} y_i y_j$. Thus one has $\sum_{i<j} y_i y_j =-1/2
\sum y_i^2$ which proves the lemma.  \findemo

Using Lemma \ref{a2}, we easily prove that there are no algebraic
curves of the same bidegrees as curves in Proposition~\ref{pseudo}
having the same positions with respect to the line pencil $\cL$.

\begin{prop}\label{alg ob}
There does not exist a nonsingular real algebraic curve $\Gamma$ of
 bidegree $(d,0)$ in $\Sigma_2$ such that there exists real numbers
 $x_1<x_2<0<z_1<x_4<z_2<x_6$ satisfying
\begin{itemize}
\item $\Gamma$ has a   maximal
tangency point with  the fibers $X-x_iZ=0$ with $i=1,\ldots 4$,
\item $\Gamma$ has $d$ distinct real intersection points with the
  fibers $X-z_iZ=0$ with $i=1,2$.\findemo
\end{itemize}  
\end{prop}
\textit{Proof. }Suppose such a real algebraic curve $\Gamma$
exists. Then, in an appropriate standard system of coordinates on
$\Sigma_2$, the curve $\Gamma$ has the following equation
$$Y^d + a_{2}(X)Y^{d-2} + a_{3}(X)Y^{d-3} + \ldots + a_{d}(X) $$ where
 $a_i(X)$ is real a polynomial of degree $2i$.  According to
 Lemma~\ref{a2}, one has $a_2(z_i)<0$, so the polynomial $a_2$ is not
 identically zero.  According to Lemma~\ref{a2}, one has $a_2(x_i)=0$,
 so the $x_i$'s are exactly the simple roots of the polynomial
 $a_2$. However, as $z_1<x_4<z_2$ and $a_2(z_1)$ and $a_2(z_2)$ have
 the same sign, there should exist an extra root $x_5$ of $a_2$ in the
 interval $]z_1;z_2[$. Hence, $a_2$ would be a non-null polynomial of
 degree 4 with at least 5 roots, which is impossible.\findemo

\section{Concluding remarks}

\begin{enumerate}
\item The patchwork in Section \ref{patchwork} can be realized
  algebraically as soon as one keeps only any three out of the four
  maximal tangency points.
Patchworking a curve with Newton polygon the triangle with vertices
$(1,0), (0,d)$ and $(2 d - 1,0)$,  one can keep two of them. Applying
  Shustin's theorems (see \cite{Sh3}, \cite{Sh1}), one can keep  the
  maximal tangency point 
  coming either from the curve $P_2$ or from the curve $P_3$. 

\item It is fairly easy to generalise the main theorem to other
rationally ruled surfaces and to construct a lot of other
examples. For instance in the same way one proves that there exists
pseudoholomorphic patchworked curves of bidegree $(d,0)$  in
$\Sigma_n$ which are not
$\cL$-isotopic to any 
such real algebraic curve as soon as   $n\ge
2$ and $d \ge 3$. One can also construct examples in $\Sigma_1$ of
bidegree $(d,0)$ for $d\ge 5$.

\item This paper is part of a work in progress in which we investigate
tangencies of curves with respect to a pencil of lines. For real
algebraic curves, one can obtain restrictions valid in any degree by
means of certain subresultants. We point out that most of these
restrictions do not hold for real pseudoholomorphic curves.  Our proof
here relies on the simplest example of algebraic prohibitions obtained
in this way.

\end{enumerate}

\small

\vspace*{2 ex}
\begin{tabular}{lll}
\textbf{Beno\^{i}t Bertrand}& \hspace{10ex}  & \textbf{Erwan Brugallé}\\
Section de math\'ematiques &    &         Université de Paris 6 Pierre
et Marie Curie\\
Universit\'e de Gen\`eve &    & 175 rue du Chevaleret \\
 case postale 64, &          &                 75 013 Paris \\
2-4 rue du Li\`evre,          &             &                France\\
Gen\`eve &&\\
Suisse &&\\
\\
E-mail : benoit.bertrand@math.unige.ch & &E-mail : brugalle@math.jussieu.fr

\end{tabular}

\end{document}